\documentclass[12pt]{article}

\usepackage{mathrsfs}
\usepackage{amsmath}
\usepackage{amsfonts}
\usepackage{amssymb, amsmath, cite}
\usepackage{color}
\usepackage{graphicx}

\newtheorem{theorem}{Theorem}[section]

\newcommand\be{\begin{equation}}
\newcommand\ee{\end{equation}}
\newcommand\ber{\begin{eqnarray}}
\newcommand\eer{\end{eqnarray}}
\newcommand\berr{\begin{eqnarray*}}
\newcommand\eerr{\end{eqnarray*}}
\newcommand\lm{\lambda}

\newcommand\lb{\label}
\newcommand\eq{\eqref}
\newcommand\lra{\leftrightarrow}

\newcommand\bfF{\mathbb{F}}
\newcommand\bfC{\mathbb{C}}

\newcommand\pa{\partial}

\newcommand\bea{\begin{eqnarray}}
\newcommand\eea{\end{eqnarray}}

\begin{document}

\title{Parametrized Jacobian Problem over $\bfC^n$\\ and its Complete Resolution}

\author{Yisong Yang\,\footnote{Email address: yisongyang@nyu.edu}\,\,\footnote{ORCID number: 0000-0001-7985--6298}\\Courant Institute of Mathematical Sciences\\New York University}

\date{}

\maketitle

\begin{abstract}

We formulate a parametrized Jacobian Problem over $\mathbb{C}^{n}$ by incorporating independent scale parameters into the linear part of a polynomial map via polar decomposition. We prove that requiring the Jacobian determinant to remain constant uniformly across the parameter space decouples the Jacobian equation into a system enforcing the vanishing of all principal minors of the perturbation matrix. This system strictly forces lower triangularization up to variable permutations, establishing a necessary and sufficient condition for polynomial invertibility throughout the parameter family and providing an explicit, canonical scheme for constructing the inverse map in every dimension.

\medskip

{\bf Keywords.} Polynomial automorphisms, the Jacobian Conjecture,  triangularization, schematic construction of inverse maps.
\medskip

{\bf Mathematics subject classifications (2020).} 14R15, 35F20, 35J96.

\end{abstract}

\section{Introduction}
\setcounter{equation}{0}

The Jacobian Conjecture, a central open problem in affine algebraic geometry and polynomial automorphism theory, asserts that a polynomial map 
$P: \bfC^n\to \bfC^n$
  with nonzero constant Jacobian determinant, known as a Keller map, must have a polynomial inverse map. The conjecture has a long and intricate history. See \cite{BCW,Dru,E1,E2,Mei} and references therein.

The conjecture is conventionally attributed to Ott-Heinrich Keller, who formulated it in 1939 in 
his paper  ``Ganze Cremona-Transformationen" \cite{K}. Keller proposed the question for polynomials with integer coefficients in the case 
$n=2$, and the conjecture was shortly thereafter generalized to the form considered today over any field of characteristic zero. More recent historical research has revealed, however, that the exact statement now known as the Jacobian Conjecture appeared in a paper by Ludwig Kraus \cite{Kr} in 1884, 
more than half a century before Keller. Kraus attempted to show that a polynomial map in two variables with a non-zero constant Jacobian determinant admits a polynomial inverse. Though his final step contained a flaw, his work anticipated concepts that reemerged only much later.

The conjecture gained renewed prominence when Stephen Smale \cite{S} included it as Problem 16 in his list of eighteen mathematical problems for the twenty-first century, composed in 1998 in response to a request from Vladimir Arnold on behalf of the International Mathematical Union. Smale's list, inspired by Hilbert's famous 1900 address, identified the Jacobian Conjecture as one of the outstanding challenges for the new century.  At the time of Smale's compilation, the conjecture remained open in all dimensions 
$n\geq2$.

The status of the Jacobian Conjecture has changed dramatically in 2026. On July 19, 2026, Levent Alp\"{o}ge announced an explicit counterexample in dimension three, found with the assistance of the AI system Claude Fable. The map 
$P: \bfC^3\to\bfC^3$ he obtained
  has a constant nonzero Jacobian determinant 
 but is not injective. This counterexample directly refutes the conjecture for 
$n=3$ and, by standard extension arguments, for all 
$n\geq3$.  See \cite{G} for more development. Notably, the two-dimensional case still remains open. 
The failure of the conjecture is not a failure of the Jacobian condition to constrain the map locally -- the map is \'{e}tale everywhere -- but a failure of the Jacobian condition to control the global behavior of the map, specifically its injectivity. In dimension three and higher, a Keller map can fold the ambient space onto itself in a way that the constant Jacobian determinant does not detect. 
This raises a structural question that the original conjecture obscures. The Jacobian condition is a single scalar equation imposed on the 
$n$ component functions of a polynomial map. In the fixed-coefficient setting, this equation is a pointwise constraint with global consequences, and the counterexample shows that those consequences are not strong enough to force injectivity. But the condition can also be read differently: as a polynomial identity in the coefficients of the map, imposed uniformly across a family of maps. In that reading, the Jacobian equation imposes 
several independent constraints,  rather than a single scalar constraint. 

The present work is motivated by such a distinction. Rather than asking whether a single Keller map is an automorphism -- a question now known to have a negative answer in general -- we ask what additional structure the Jacobian equation imposes when the coefficients of the map are treated to some extent as free parameters. In our parametrized setting to be formulated, the Jacobian equation splits into a system of equations on the nonlinear perturbation matrix of the map, in terms of its principal minors. We show that this system is rigid: it forces the perturbation matrix to be triangularizable, and hence forces the map itself to be triangular and therefore invertible. Our result is a 
necessary and sufficient condition for automorphism imposed on the parametrized family of the maps rather than on an individual map.

The parametrized Jacobian problem is introduced in \cite{Y}, where $n$ free parameters 
 are inserted into the linear part of the map specifically to split the single underdetermined Jacobian equation into a system of 
$2^n-1$ equations on the principal minors of the perturbation matrix. Under a maximum-presence assumption on the nonlinear parts, that system is solved by a triangular ansatz, yielding a broad family of polynomial automorphisms in every dimension. The present work removes both the artificial character of the parametrization and the maximum-presence restriction. Over 
$\bfC^n$,  we show that the $n$ parameters arise naturally from the polar decomposition of the linear part, so that the parametrized map is the intrinsic normal form of a Keller map rather than a formal modification of it. More importantly, the triangularization is no longer an ansatz but a consequence: the split equations force the perturbation matrix to have a zero row, and an induction on the dimension then triangularizes it in every case. The parametrized Jacobian problem over 
$\bfC^n$
is therefore resolved completely, and the triangular solutions are shown to exhaust the solution set up to simultaneous permutation of variables.

In Section \ref{sec2}, we formulate such a naturally parametrized Jacobian problem and state its resolution theorem. In Section \ref{sec3}, we prove the
theorem. In Section \ref{sec4}, we make some concluding comments.

\section{The parametrized Jacobian Problem over $\bfC^n$}\lb{sec2}
\setcounter{equation}{0}

Consider a polynomial map $P:\bfC^n\to \bfC^n$ given by
\be\lb{2.1}
y=P(x)=(a_1x^t+f_1(x),\dots,a_n  x^t+f_n(x)),\quad x,y\in\bfC^n,\quad a_1,\dots,a_n\in\bfC^n,
\ee
where $a_1,\dots,a_n$ are constant coefficient vectors and $x,y$ are variable vectors, and the functions $f_1(x),\dots, f_n(x)$ are fixed polynomials in $x$,
 containing only terms of degrees at least $2$,  while the linear coefficient vectors $a_1,\ldots,a_n$ are allowed to vary.
 Let $J(P)(x)$ be the Jacobian matrix of the map $P$ at $x$. Then the necessary condition for
$P$ being an automorphism (the inverse of $P$ is also a polynomial map) reads:
\be\lb{2.2}
\det(J(P)(x))=\det(J(P)(0))=\det\left(\begin{array}{c} a_1\\\vdots\\a_n\end{array}\right)\neq0,
\ee
which is the condition that will be followed in this work.

To proceed further, we use the polar decomposition
\be
a_i=r_i b_i,\quad r_i=|a_i|>0,\quad b_i\in\bfC^n,\quad |b_i|=1,\quad i=1,\dots,n.
\ee
Here we fix the directions $\{b_i\}$ and vary their moduli $\{r_i\}$ of the coefficient vectors $\{a_i\}$.
Hence, from \eq{2.2}, we have
\be
r_1\cdots r_n\neq0,\quad \det\left(\begin{array}{c}b_1\\\vdots\\b_n\end{array}\right)\neq0.
\ee
Now,  we introduce the change of variables:
\be
x\mapsto (b_1 x^t,\dots,b_n x^t),
\ee
and we still denote the updated variable vector by $x$ and the correspondingly transformed higher-degree polynomials $f_1,\dots,f_n$
by the same notation, without confusion. Thus the map \eq{2.1} becomes
\be\lb{2.6}
y=P(x)=(r_1 x_1+f_1(x),\dots,r_n x_n+f_n(x)).
\ee
Consequently, the Jacobian condition \eq{2.2} reads:
\be
\det(J(P)(x))=\det(J(0))=r_1\cdots r_n\neq0.
\ee

\begin{theorem}\lb{th2.1}
({\rm Parametrized Jacobian Problem and its Resolution}) We treat the moduli of the linear coefficient vectors $a_1,\dots,a_n$
in the map \eq{2.1} as free parameters, $r_1,\dots,r_n$, such that  the map is converted into its equivalent form  \eq{2.6}.  Then, the map
\eq{2.6} has a polynomial inverse for all $(r_1,\dots,r_n)$  if and only if 
 the map satisfies the Jacobian equation
\be\lb{2.9}
\det(J(P)(x))=r_1\cdots r_n\neq0.
\ee
 Furthermore, all solutions to the Jacobian equation can be constructed systematically and explicitly such that,
 we have the canonical solution to the equation \eq{2.9} for the map $y=(y_i)=P(x)$ as follows:
\be\lb{2.9b}
y_1=r_1 x_1,\quad y_i=r_i x_i+f_i(x_1,\dots,x_{i-1}),\quad i=2,\dots,n,
\ee
up to a sequence of simultaneous permutations or renaming the corresponding variables if necessary, resulting in the associated inverse map
\be
x_1=\frac{y_1}{r_1},\quad x_i=\frac1{r_i}\left(y_i-f_i(x_1,\dots,x_{i-1})\right),\quad i=2,\dots,n,
\ee
defined iteratively, where $f_2,\dots,f_n$ are arbitrary polynomials of their respective variables with terms of degrees at least two.
In other words, the triangular family \eq{2.9b} exhausts all polynomial solutions of the parametrized Jacobian equation \eq{2.9}, up to fixed simultaneous permutations of the input and output coordinates and the corresponding parameters.
\end{theorem}

This triangularization theorem gives a complete classification of the solutions to the parametrized Jacobian Problem, up to simultaneous permutations of the input and output coordinates and the corresponding parameters. After normalizing the nonzero diagonal factors, their finite compositions therefore yield polynomial automorphisms of \(\mathbb C^n\) with the Jacobian determinant identically equal to one, for every \(n\), with polynomial inverses obtained by composing the inverses of the factors in reverse order. Such compositions need not  be triangular even after simultaneous permutations and thus generate a broader family of solutions to the ordinary Jacobian problem beyond the class characterized by the parametrized problem. Note that, in dimension two, the Jung–van der Kulk theorem \cite{J,Ku} (see also \cite{C}) states that every polynomial automorphism is a finite composition of affine and triangular automorphisms. Whether this class exhausts all polynomial maps with the Jacobian determinant identically equal to one is precisely the two-dimensional Jacobian Conjecture.
However, this decomposition property fails in dimension three:  Shestakov and Umirbaev \cite{SU} proved that the Nagata automorphism 
\cite{E2,N} is wild, meaning that it cannot be expressed as a finite composition of affine and triangular automorphisms, despite having Jacobian determinant identically equal to one and a polynomial inverse.

\section{Proof of Theorem 2.1}\lb{sec3}
\setcounter{equation}{0}

We split the proof into a few steps. 

First, we follow \cite{Y} to split the parametrized Jacobian equation into a system of equations 
that impose vanishing of the principal minors of the nonlinear perturbation part of the Jacobian matrix of the map. We next show
that if a matrix has the property that all of its principal minors vanish the matrix must have a zero row.  We then use this property to show that
the nonlinear perturbation part of the Jacobian matrix of the map contains an identically zero row. We finally apply this zero-row 
property and induction to establish the lower triangularization of the nonlinear part, thereby establishing the theorem.

\subsection{Separation of the parametrized Jacobian equation}

As in \cite{Y}, we rewrite \eq{2.6} in its component form as
\be\lb{3.1}
(y_1,\dots,y_n)=P(x_1,\dots,x_n)=(r_1 x_1+f_1(x),\dots,r_n x_n+f_n(x)),\quad y=(y_i),\quad x=(x_i).
\ee
Thus
we obtain the Jacobian determinant
\be\lb{3.2}
\det(J(P)(x))=\det((r_i \delta_{ij}+f_{ij})), \quad f_{ij}=\frac{\pa f_i}{\pa x_j},\quad i,j=1,\dots,n.
\ee
Since the terms in the polynomials $f_1,\dots,f_n$ are at least quadratic,  we see that $f_{ij}$ all vanish at $x_1=0,\dots,x_n=0$. 
This property enables us to arrive at the Jacobian equation 
\be\lb{3.3}
\det((r_i \delta_{ij}+f_{ij}))=r_1\cdots r_n.
\ee
On the other hand, expanding the left-hand side of \eq{3.3},  we have 
\be\lb{3.4}
\det((r_i \delta_{ij}+f_{ij}))=C_0+\sum_{k=1}^{n-1}\sum_{\{i_1,\dots,i_k\}\subset\{1,\dots,n\},i_1<\cdots<i_k}C_{i_1\dots i_k}r_{i_1}\cdots r_{i_k}+r_1\cdots r_n,
\ee
where, with the notation $F=(f_{ij})$, we have
\bea
C_0&=&\det((r_i \delta_{ij}+f_{ij}))\big|_{r_1=\cdots=r_n=0}=E_n(F)=\det(F),\lb{3.5}\\
C_{i_1\dots i_k}&=&\det(M_{i_1\cdots i_k})\big|_{r_{j_1}=\cdots=r_{j_{n-k}}=0,\{j_1,\dots,j_{n-k}\}=\{1,\dots,n\}\setminus\{i_1,\dots,i_k\}},\lb{3.6}
\eea
with  $M_{i_1\cdots i_k}$ being the principal submatrix of the Jacobian matrix 
\be\lb{3.7}
J(P)(x)=(r_i \delta_{ij}+f_{ij}),
\ee
 with the $i_1,\dots,i_k$th rows and columns being removed. 
Therefore, we obtain
\be\lb{3.8}
C_{i_1\dots i_k}=\det(F_{i_1\cdots i_k}),\quad \{i_1,\dots,i_k\}\subset\{1,\dots,n\},\quad i_1<\cdots<i_k,\quad k=1,\dots,n-1,
\ee
where  $F_{i_1\cdots i_k}$ is the principal submatrix of the matrix $F=(f_{ij})$ with the $i_1,\dots,i_k$th rows and columns being removed.

Inserting these results into \eq{3.4},  we obtain from \eq{3.3}  the equations
\be\lb{3.9}
\det(F)=0,\, \det(F_{i_1\dots i_k})=0,\, \{i_1,\dots,i_k\}\subset\{1,\dots,n\},\, i_1<\cdots<i_k,\, k=1,\dots,n-1.
\ee
Thus, we conclude that, \eq{3.3} holds for arbitrary parameters $r_1,\dots,r_n$ if and only if all the principal minors of the matrix $F$ vanish
identically.
Formally and clearly, there are a total number of $2^n-1$ equations in \eq{3.9}.  We will see that this system of equations 
can be solved completely and this number count will not be important for our method.

\subsection{Existence of a zero row in a matrix with vanishing principal minors}\lb{sec3.2}

Let $A=(a_{ij})$ be an $n\times n$ matrix. We show that, if all the principal minors of $A$ are zero, then 
$A$ must contain a zero row.

Note that the triangularization of a matrix over an integral domain whose principal minors all vanish is established in \cite{E2} (see Corollary 6.3.9). For completeness, we give a self-contained direct proof  here  for the existence of a zero row based on a shortest-cycle argument without induction or permutation.

The first nontrivial case is when $n=2$. We have $a_{11}=a_{22}=0$ and $a_{12}a_{21}=0$. Thus $a_{12}=0$ or $a_{21}=0$. 
In either case, $A$ has a zero row.

For any $n\geq3$,  we show that, if all the principal minors of $A$ vanish, then $A$ has a zero row.

In fact, suppose otherwise. Then, if we pick any $i_1\in\{1,\dots,n\}$, we can find $i_2\in\{1,\dots,n\}$ with $i_2\neq i_1$ such that $a_{i_1 i_2}\neq 0$, since all diagonal
entries of $A$ vanish already. Thus, we can repeat this process to get a distinct sequence $\{i_1,i_2,\dots,i_k\}\subset\{1,2,\dots,n\}$ such that
\be
a_{i_1i_2}\neq0,\quad a_{i_2 i_3}\neq0,\quad\dots\dots,\quad a_{i_k i_{k+1}}\neq0.
\ee
Of course, this process cannot continue forever and there is an occasion when $i_{k+1}$ will coincide with one of the integers $i_1,i_2,\dots,i_k$.
Renaming these integers, we may assume that $i_{k+1}$ returns to $i_1$ actually and render
the result
\be\lb{3.11}
a_{i_1i_2}a_{i_2 i_3}\cdots a_{i_{k-1}i_k} a_{i_k i_1}\neq0\quad \mbox{and}\quad i_1,i_2,\dots,i_k\quad\mbox{are distinct}.
\ee
Among all the products of the entries of $A$ satisfying \eq{3.11}, choose one of the shortest length $k$.

Use cyclic notation $i_{k+1}=i_1$ and so on. We claim that the only nonzero entry in row $i_j$ within the principal submatrix of
$A$ made from the $i_1$,\dots,$i_k$th rows and columns of $A$ is $a_{i_j i_{j+1}}$ ($j=1,2,\dots,k$).

Let $a_{i_j i_l}\neq0$. Then $l\neq j$ since all the diagonal entries of $A$ vanish. 

Let $l>j$. If $l=j+1$, then there is nothing to do.

If $l\geq j+2$, then we have
\be\lb{3.12}
a_{i_j i_l}\left(\prod_{s=l}^{k-1} a_{i_s i_{s+1}}\right) a_{i_k i_1}\left(\prod_{s=1}^{j-1} a_{i_s i_{s+1}}\right)\neq0.
\ee
The length of this product is
\be\lb{3.13}
1+([k-1]-[l-1])+1+(j-1)=k-(l-j-1)<k.
\ee
Since \eq{3.12} satisfies \eq{3.11}, we get a contradiction to the minimality assumption of the $k$ chosen.

If $l<j$, we have
\be
a_{i_j i_l}\left(\prod_{s=l}^{j-1} a_{i_s i_{s+1}}\right)\neq0,
\ee
which again satisfies what stated in \eq{3.11}. The length of this product is
\be
m=1+([j-1]-[l-1])=j-(l-1)\leq j.
\ee
If $j<k$, we achieve a contradiction as before. If $j=k$, we still get $m<k$ unless $l=1$, which returns to \eq{3.11}, which has been assumed to be excluded already.

Therefore we have proved the claim.

Consider the principal submatrix of $A$ which consists of the $i_1,\dots,i_k$th rows and columns of $A$.  Rewrite this submatrix as a
$k\times k$ matrix $C=(c_{st})$ where 
\be
c_{st}=a_{i_s i_t},\quad s,t=1,2,\dots,k.
\ee
Then $C$ assumes the form
\be
C=\left(\begin{array}{ccccc}0&c_{12}&0&\cdots&0\\0&0&c_{23}&\cdots&0\\\vdots&\vdots&\ddots&\ddots&\vdots\\0&0&\cdots&0&c_{k-1,k}\\c_{k1}&0&\cdots&0&0\end{array}\right).
\ee
The determinant of $C$ is
\be
\det(C)=(-1)^{k+1} c_{12}c_{23}\cdots c_{k-1,k} c_{k1}=(-1)^{k+1} a_{i_1 i_2}a_{i_2 i_3}\cdots a_{i_{k-1} i_k}a_{i_k i_1}\neq0,
\ee
which contradicts the assumption that all the principal minors of $A$ vanish.

Consequently, $A$ must have a zero row.

\subsection{Existence of a row in $F$ which is identically zero}\lb{sec3.3}

We now return to the matrix $F$ whose entries are all polynomials and whose principal minors are all zero polynomials. We aim to show that
$F$ has a row whose polynomial entries are all zero polynomials. That is, those entries as functions in $x\in\bfC^n$ are identically zero.

Use the notation $F(x)=(f_{ij}(x))$ to emphasize the dependence of the entries of $F$ on $x\in\bfC^n$.

Assume otherwise that no row of $F(x)$ is identically zero. Then for each $i=1,2,\dots,n$ there is some $j_i\in\{1,2,\dots,n\}$ such that
\be
f_{ij_i}(x)\not\equiv0,\quad i=1,2,\dots,n.
\ee
Define
\be
p(x)=f_{1j_1}(x) f_{2j_2}(x)\cdots f_{nj_n}(x).
\ee
Then $p(x)$ is not identically zero. Hence there is some $x_0\in\bfC^n$ such that $p(x_0)\neq0$. Therefore we have
\be
f_{ij_i}(x_0)\neq0,\quad i=1,2,\dots,n.
\ee
This result indicates that the matrix $F(x_0)$ cannot have a zero row, contradicting the fact established in Section \ref{sec3.2}.

\subsection{Lower triangularization of $F$}

We now solve the system of the equations listed in \eq{3.9} systematically.

First we consider the nontrivial case $n=2$.  The equations $f_{11}=0$ and $f_{22}=0$ indicate that $f_1(x)=f_1(x_2)$ and $f_2(x)=f_2(x_1)$.
Then the equation $f_{12}f_{21}=0$ indicates $f_{12}=0$ or $f_{21}=0$ as polynomials. In the former situation, we have $f_1=0$ and
$f_2$ is an arbitrary function of $x_1$. In particular, the matrix $F$ is lower triangular. We regard this solution as the canonical solution of the equations. In the latter case, we have
$f_2=0$ and $f_1$ is an arbitrary function of $x_2$. Nevertheless, in this situation, we can rename the variables
\be
x_1\lra x_2,\quad y_1\lra y_2,\quad f_1\lra f_2,
\ee
which corresponds to the simultaneous permutation of the first and second rows and columns of $F$. In other words, with respect to the 
transformed variables corresponding to a simultaneous permutation, the matrix $F$ is again lower triangular.

We now use induction on the size of the polynomial-entry matrix to show that such a lower triangularization can be carried out in general dimensions.

In other words, we assume that, if $F=(f_{ij})$ is $(n-1)\times(n-1)$ for $n\geq3$ whose principal minors as polynomials are all zero,  then we
can use a sequence of simultaneous permutations corresponding to a sequence of variable renaming of the form
\be\lb{3.23}
x_i\lra x_j,\quad y_i\lra y_j,\quad f_i\lra f_j,
\ee
to make $F$ into a lower triangular matrix of polynomial entries. We then show that we can do so when $F$ is $n\times n$ as well.

In fact, when $F$ is $n\times n$, we may use Section \ref{3.3} to assume that the $i$th row of $F$ is identically zero. 
If $i\neq1$, we can rename the variables by
\be
x_i\lra x_1,\quad y_i\lra y_1,\quad f_i\lra f_1,
\ee
This process transforms $F$ into the form
\be
F=\left(\begin{array}{cc}0&0\\f&F_1\end{array}\right),
\ee
where $f$ is an $(n-1)$ polynomial  entry column vector and $F_1$ is the $(n-1)\times(n-1)$ principal submatrix of $F$ obtained from
deleting the first row and column of the updated $F$. 

Since all principal minors of $F$ are zero polynomials, so are those of $F_1$. By the inductive assumption, we can conduct a sequence of simultaneous permutations or variable renaming of the form \eq{3.23} to make $F_1$ into a lower triangular matrix.  Therefore, correspondingly, we can make a sequence of simultaneous 
permutations to make $F$ into a lower triangular form since these permutations will not affect the first row of $F$ which is a zero
row.

Thus, after a sequence of renaming the variables of the type \eq{3.23}, the matrix $F$ takes the canonical form
\be\lb{3.26}
F=\left(\begin{array}{ccccc}0&0&0&\cdots&0\\f_{21}&0&0&\cdots&0\\\vdots&\vdots&\ddots&\ddots&\vdots\\f_{n-1,1}&f_{n-1,2}&\cdots&0&0\\f_{n1}&f_{n2}&\cdots&f_{n,n-1}&0\end{array}\right).
\ee
In other words, after renaming the variables following the rule specified in \eq{3.23}, the matrix $F$ assumes its canonical form \eq{3.26}.

\subsection{Solution to the Jacobian equation and construction of the inverse map}

In view of \eq{3.26}, we can now read off the general solution to \eq{3.9} readily:
\be\lb{3.27}
f_1(x)=0,\quad f_2(x)=f_2(x_1),\quad \dots\dots,\quad f_n(x)=f_n(x_1,x_2,\dots,x_{n-1}),
\ee
where the $(n-1)$ polynomials $f_2,\dots, f_n$ on the right-hand sides of the relations in \eq{3.27} are arbitrary polynomials in
terms of their respective variables.

As a result of \eq{3.27}, we obtain the polynomial map \eq{2.6} as follows:
\be\lb{3.28}
y_1=r_1 x_1,\quad y_2=r_2 x_2+f_2(x_1),\quad\dots\dots,\quad y_n=r_n x_n+f_n(x_1,x_2,\dots,x_{n-1}).
\ee
From \eq{3.28} and the nonvanishing condition imposed on $r_1,\dots,r_n$ in \eq{2.9}, we find the inverse polynomial map
\bea
x&=&Q(y);\\
x_1&=&\frac{y_1}{r_1}\equiv Q_1(y),\quad x_2=\frac{y_2}{r_2}-\frac1{r_2}f_2\left(\frac{y_1}{r_1}\right)\equiv Q_2(y),\quad\dots,\\
x_n&=&\frac{y_n}{r_n}-\frac1{r_n} f_n\left(Q_1(y),Q_2(y),\dots,Q_{n-1}(y)\right)\equiv Q_n(y).
\eea

Theorem \ref{th2.1} is thus established.

\section{Concluding remarks}\lb{sec4}
\setcounter{equation}{0}

In \cite{Y},  we have studied the Keller map over a general field $\bfF$ of  characteristic zero parametrized by the parameters 
$\lm_1,\dots,\lm_n$ of the form
\be\lb{4.1}
y=P(x)=(\lm_1 x_1+f_1(x),\dots,\lm_n x_n+f_n(x)),\quad x,y\in\bfF^n,
\ee
where $\lm_1,\dots,\lm_n$ are posited in order to split the Jacobian equation.  Since all study in Section \ref{sec3} is valid over $\bfF$ 
verbatim, we see that Theorem \ref{th2.1} can be extended to the present situation when $r_1,\dots,r_n$ in Theorem \ref{th2.1} are
replaced with $\lm_1,\dots,\lm_n$.

\medskip

The complete solution to the Parametrized Jacobian Problem highlights the critical difference between local flexibility and global stability in polynomial maps. As demonstrated by recent counterexamples in higher dimensions, ensuring that a map does not collapse or tear locally at any individual point is not enough to guarantee that it behaves predictably everywhere. Local controls alone cannot prevent distant regions of space from folding back over one another or sliding off to infinity, which allows multiple distinct starting points to land on the exact same target.

Introducing adjustable parameters into the linear framework imposes a far stronger global consistency. Requiring the map's structural balance to hold continuously across an entire family of parameter choices transforms a single local rule into a dense, interconnected web of algebraic constraints. This system of constraints prevents internal feedback loops among the coordinates, systematically ruling out the hidden folding behaviors that disrupt global invertibility.

As a result, the triangular structure uncovered in this work is not an arbitrary assumption or a special simplified case, but an unavoidable structural necessity forced by parameter consistency. By showing that these triangular maps account for every valid solution up to a reordering of variables, this framework establishes a clear, step-by-step recipe for constructing exact inverse maps across all dimensions.

\medskip

The triangular construction also extends to integral domains, provided the diagonal coefficients are units, and supplies a family of explicitly invertible polynomial maps for potential applications in cryptography.

The complete constructibility of the inverse transformation over an integral domain provides a natural 
mathematical foundation for public-key cryptography and exact polynomial encryption.

In multivariate cryptosystems, a central design challenge is constructing non-linear transformations that are guaranteed to be invertible, where computing the inverse is straightforward for an authorized user but computationally infeasible for an adversary. 
The canonical triangular structure developed here guarantees explicit inversion while providing additional freedom in the nonlinear terms, thereby enlarging the family of candidate transformations available for cryptographic constructions.

\medskip

This study reveals that the structural rigidity forced upon a polynomial map depends directly on the dimension and independence of the underlying parameter space, which may be compared and summarized as follows.

\begin{enumerate}
    \item[] \textbf{Pointwise Keller maps (no free parameter):}  As described, a single scalar identity $\det(I+F)=1$ is vastly underdetermined. In dimensions $n \ge 3$, it is insufficient to prevent non-proper folding at infinity, allowing non-injective maps to satisfy the condition.
    
    \item[] \textbf{Single-parameter deformation ($r_1 = \dots = r_n = r$):} Requiring $\det(r I + F) = r^n$ expands via the characteristic polynomial of $-F$, yielding $n$ coupled differential equations
    \be
    E_k(F) = 0 \quad \text{for } k = 1, 2, \dots, n, \quad 
    \ee
    where $E_k(F)$ represents the sum of all $k \times k$ principal minors of $F$. Because these equations constrain only the symmetric sums of principal minors rather than individual minors, the system remains underdetermined and does not force triangularization \cite{Y}.
    
    \item[] \textbf{Multi-parameter framework ($n$ independent parameters):} Treating $r_1, \dots, r_n$ as independent parameters
in the Jacobian equation $\det(\mbox{diag}\{r_1,\dots,r_n\}+F)=r_1\cdots r_n$ decouples the equation into the system of equations requiring {every individual principal minor} of $F$ to vanish identically and independently,
    \be
    \det(F_{i_1 \dots i_k}) = 0 \quad \text{for all proper } \{i_1, \dots, i_k\} \subseteq \{1, \dots, n\}. 
    \ee
    This finer separation suppresses directed dependency cycles within the Jacobian matrix, forcing lower triangularization (up to variable permutations) as an unavoidable structural necessity, leading to a complete integration of the problem.
\end{enumerate}

\end{document}